\documentclass[12pt]{article}
\usepackage{amsfonts}

\newtheorem{cor}{Corollary}
\newtheorem{lm}{Lemma}
\newtheorem{thm}{Theorem}

\newcommand{\ov}{\overline }

\begin{document}

\title{On arithmetical and dynamical properties of Lorentz maps of the torus}
\author{S.~Aranson\and E.~Zhuzhoma}
\date{}
\maketitle

\begin{abstract}
We proved that any Lorentz transformation of 2-torus $T^2$ is Anosov automorphism. One completely describes admissible parameters of Lorentz transformations and their arithmetical properties. One proved that an admissible speed light parameter has a countable spectra accumulating to this parameter.
\end{abstract}

\section{Introduction}

Lorentz manifold $M$ is a pseudo-Riemannian manifold with the signature $(1, n-1)$, where $n = \dim M$. This means that Lorentzian metric is defined in each fiber $T_pM$, $p\in M$, of the tangent bundle by the nondegenerate bilinear form $T_pM\times T_pM\to \mathbb R$ with the diagonal type $(-, +\ldots +)$. In contrast with a Riemannian metric which can be defined on any smooth compact manifold, a Lorentzian metric is defined with certain topological conditions. For example, a compact manifold endowed with a Lorentzian metric must have a zero Euler characteristic \cite{BeemEhrlich-book81}. Isometry groups of compact Lorentzian manifolds are considered in \cite{AdamsStuch97}, \cite{DAmbra88} ,where there are many references on other articles. 2-dimensional torus is a unique closed orientable surface admitting the structure of Lorentzian manifold. Moreover, the torus admits the Lorentzian metric of constant zero curvature. It is well known that any Lorentzian manifold with a metric of constant zero curvature is locally isomorphic to the Minkowski space-time \cite{Wolf-book72}. Therefore, it is natural to consider the isometry group of the 2-torus that is endowed by the Minkowski metric. The group isometry of such Lorentzian torus is richer, in sense, than the isometry group of Riemannian torus. The later consists of the rigid Euclidean shifts and the isometry exchanging ``meridian" and ``parallel".  In addition with rigid Euclidean shifts, the Lorentzian torus contains the 2-parametric family of Lorentz transformations. Lorentz maps play an important role in electrodynamics, Special Theory of Relativity, and other physic branches, and they have the interest itself.

In this paper, we study dynamical properties of Lorentz transformations on the 2-torus. We completely describe the  admissible parameters of Lorentz transformations and their arithmetical properties. At the end of the paper we present the algorithm for calculating a spectra of fixed admissible speed light parameter.

{\it Acknowledgment}. The research was partially supported by RFFI grant 02-01-00098. We thanks Ludmila Litinskaya and Igor Aranson for valuable discusions. This work was done while the second author was visiting Rennes 1 University (IRMAR) in March-Mai 2004. He thanks the support CNRS which made this visit possible. He would like to thank Anton Zorich and Vadim Kaimanovich for their hospitality.

\section{Statement of main results}

Let us represent the 2-torus $T^2$ as the factor-space $\mathbb{R}^2/\mathbb{Z}^2$ where $\mathbb{R}^2$ is Euclidean plane endowed with coordinates $(x,t)$, and $\mathbb{Z}^2$ is the group of covering transformations that consists of integer shifts. The natural projection $\pi : \mathbb{R}^2\to \mathbb{R}^2/\mathbb{Z}^2$ is a universal cover. This covering map and usual Euclidean metric on $\mathbb{R}^2$ induces on $T^2$ the structure of a flat Riemannian manifold. Denote such torus by $T^2_E$. The Minkowski metric generated by the form 
$ds^2 = -dx^2 + c^2dt^2$ induces on $T^2$ the metric of Lorentzian manifold of constant zero curvature. Denote such torus by $T^2_L$.

A homeomorphism $T^2_L\to T^2_L$ is called {\it Lorentz transformation} if its covering map is of the following type
\begin{equation}\label{lorentz-map}
\overline x = \frac{x-Vt}{\sqrt{1-\frac{V^2}{c^2}}},\quad 
\overline t = \frac{t-\frac{V}{c^2}x}{\sqrt{1-\frac{V^2}{c^2}}},\quad {\mbox where}\quad 0<\vert V\vert < c.
\end{equation}
Such a homeomorphism is denoted by $\phi_{V,c}$. The parameters $(V,c)$ for which there exist a Lorentz transformation $\phi_{V,c}$ are called {\it admissible}. Sometimes we'll say that the pair $(V,c)$ is admissible. It is easy to see that the set of admissible pairs $(V,c)$ is not empty. One can check that a Lorentzian transformation is an isometry of $T^2_L$. 

The following theorem describes dynamical properties of Lorentz transformations on the both $T^2_E$ and $T^2_L$.
\begin{thm}\label{dynamic-properties}
Every Lorentz transformation $\phi_{V,c}$ is  an Anosov hyperbolic automorphism of $T^2_E$. Any point with rational coordinates is periodic. Other points are dense under $\phi_{V,c}$ on the both $T^2_E$ and $T^2_L$.
\end{thm}

Note that if the pair $(V,c)$ is admissible, then the pair $(-V,c)$ is also admissible. For every fixed admissible $c$, the set of values $V$ that form the admissible pairs $(V,c)$ is called a {\it spectra} of $c$. The following theorem describes arithmetical properties of admissible parameters and specters.
\begin{thm}\label{arithmetical-properties}
A positive value is an admissible (light speed) $c$ if and only if this value is an irrational algebraic number which is the square root from a rational number (thus,  on the positive halfline, there are the dense subset of admissible $c$). Moreover,
\begin{enumerate}
\item A positive value is an admissible (relative speed) $V$ if and only if this value is a rational of the type $\frac{n}{m}$ such that $p = \frac{m^2-1}{n}$ is an integer, where $m$, $n\in \mathbb{N}$ are coprime.
\item For every admissible $V$, there is a unique admissible $c\stackrel{\rm def}{=}c(V)$ such that they form the admissible pair $(V, c(V))$.
\item  Given any admissible $c$, there is a countable spectra 
$$\ldots , -V_i, \ldots , -V_1, V_1, \ldots , V_i, \ldots $$
with $c = c(V_i)$ for every $i\in \mathbb{N}$.
\item Each spectra $\ldots , c(-V_i), \ldots , c(-V_1), c(V_1), \ldots , c(V_i), \ldots $ has exactly two points $-c(V_i)$, $c(V_i)$ of accumulation.
\item Specters that correspond to different light speeds are mutually disjoint.
\end{enumerate}
\end{thm}
As a consequence of items (3) and (4) of theorem \ref{arithmetical-properties}, we get
\begin{cor}
The set of admissible $V$ is dense on the real line.
\end{cor}

\section{Proof of main theorems}

First of all, let us consider when a homeomorphism $\ov \phi_{V,c}: \mathbb{R}^2\to \mathbb{R}^2$ is a covering map for some homeomorphism $\phi_{V,c}$ of the torus. Due to (\ref{lorentz-map}), the matrix
\begin{displaymath}
\mathbf{A} =
\left(
\begin{array}{cc}
a & -aV \\
-\frac{aV}{c^2} & a
\end{array}
\right)
\end{displaymath}
must be unimodal and integer, where $a = \frac{1}{\sqrt{1-\frac{V^2}{c^2}}}$. The right calculation gives that $\det A = 1$. Then the unique condition is that the elements of $\mathbf A$ are integers. Denote
\begin{equation}\label{redenote}
a\stackrel{\rm def}{=}m,\quad -aV^2\stackrel{\rm def}{=}-n,\quad -\frac{aV^2}{c^2}\stackrel{\rm def}{=}-p.
\end{equation}
It follows from $0 < \vert V\vert < c$ that $m\ge 2$. Since a pair $(V,c)$ is admissible if and only if the pair $(-V,c)$ is also admissible, later on, we assume $V > 0$. Hence, $m$, $n$, and $p$ are natural number that satisfies the following equation
\begin{equation}\label{Pell-equation}
m^2 - np = 1\quad m\ge 2, \quad n, p\in \mathbb{N}.
\end{equation}
If we know the solution $(m,n,p)$ of (\ref{Pell-equation}), then one can find the parameters $c$ and $V$:
\begin{equation}\label{cV-from-mnp}
c = \sqrt{\frac{n}{p}},\quad V = \frac{n}{m},
\end{equation}
and vise versa, given parameters $c$ and $V$, one can define the solution $(m,n,p)$ of (\ref{Pell-equation}):
\begin{equation}\label{mnp-from-cV}
m = \frac{1}{\sqrt{1-\frac{V^2}{c^2}}},\quad n = \frac{V}{\sqrt{1-\frac{V^2}{c^2}}},\quad 
p = \frac{V}{c^2\sqrt{1-\frac{V^2}{c^2}}}.
\end{equation} 
Thus, the existence of admissible parameters reduses to the existence of integer positive solutions of equation (\ref{Pell-equation}).

\begin{lm}\label{solutions-exist}
There are infinitely many admissible pairs $(c_m,V_m)$, $(c_{-m},V_{-m})$ such that
$$ c_m, V_m \to +\infty ,\quad c_{-m}, V_{-m} \to 0\quad \mbox{as}\quad m\to +\infty .$$
\end{lm}
{\it Proof}. Equation (\ref{Pell-equation}) has infinitely many integer solutions $(m,n,1)$, where $m\ge 2$, 
$n = m^2 - 1$. It follows from (\ref{cV-from-mnp}) and (\ref{Pell-equation}) that the corresponding admissible parameters $c_m = \sqrt{m^2-1} $ and $V_m = \frac{m^2-1}{m}$ tend to $+\infty$ as $ m\to +\infty $. Similarly, (\ref{Pell-equation}) has infinitely many integer solutions $(m,1,p)$, where $m\ge 2$, $p = m^2 - 1$. The corresponding admissible parameters $ c_{-m} = \frac{1}{m^2-1} $ and $V_{-m} = \frac{1}{m}$ tend to $0$ as 
$ m\to +\infty $. This proves the lemma. $\Box$

As a consequence, there are infinitely many admissible parameters $c$ and $V$.
\begin{lm}\label{c-is-irrational}
Any admissible $c$ is algebraic irrational, and any admissible $V$ is rational.
\end{lm}
{\it Proof}. It follows from (\ref{Pell-equation}) and (\ref{cV-from-mnp}) that $c = \frac{\sqrt{m^2-1}}{p}$. This implies that $c$ is algebraic irrational. The rationality of $V$ follows from (\ref{cV-from-mnp}). $\Box$
\medskip

{\it Proof of theorem \ref{dynamic-properties}.} The straight calculation gives that the eigenvalues of the matrix $\mathbf A$ are
$$ \lambda _1 = \frac{\sqrt{1-\beta ^2}}{1+\beta}\in (0;1),\quad 
\lambda _2 = \frac{\sqrt{1-\beta ^2}}{1-\beta}\in (1;+\infty ),$$
where $ \beta = \frac{V}{c} $. The corresponding eigenvectors are
$$ \vec e_1 = (c\mu , \mu ),\quad \vec e_2 = (-c\mu , \mu ),$$
where $\mu = const$. By lemma \ref{c-is-irrational}, $c$ is irrational. This implies that the corresponding transformation $\phi_{V,c}$ is  an Anosov hyperbolic automorphism of $T^2_E$. Since $\phi_{V,c}$ is linear, any point with rational coordinates is periodic, and other points are dense under $\phi_{V,c}$ on the both $T^2_E$ and $T^2_L$ \cite{Robinson-book99}. The theorem is proved. $\Box$.

In order to proof theorem \ref{arithmetical-properties}, we consider properties of solutions $(m,n,p)$ of equation (\ref{Pell-equation}) in details.
\begin{lm}\label{coprime}
Let $(m,n,p)$ is a solution of equation (\ref{Pell-equation}). Then every pair $(m,n)$, $(m,p)$, and $(m,np)$ are coprime numbers.
\end{lm}
{\it Proof}. First, let us prove that the numbers $m$, $n$ are coprime. Suppose the contradiction. Then there is the integer $k\neq 0$ such that $m = km_0$, $ n = kn_0 $. Deviding the quality $m^2 - np = 1$ on $k$, we get
$$ m\frac{m}{k} - p\frac{n}{k} = \frac{1}{k} $$
or $ m_0m - n_0p = \frac{1}{k}$. Since the left side is an integer and the right side is fractional, we get the contradiction.

The proof is similar for co-primeness of $m$ and $p$. It follows from the quality $m^2 - np = 1$ that $m$ and $np$ are also coprime. $\Box$

Admissible parameters $c$ and $V$ are said to be corresponding to the solution $(m,n,p)$ of equation (\ref{Pell-equation}) if the relations (\ref{cV-from-mnp}), (\ref{mnp-from-cV}) hold.
\begin{lm}\label{absolutly-different-for-c}
Let $(m,n,p)$, $(m^{\prime},n^{\prime},p^{\prime})$ are different solutions of equation (\ref{Pell-equation}) corresponding to the same admissible $c$ (that is, the solutions don't coinside in at least one position). Then 
$m\neq m^{\prime}$, $n\neq n^{\prime}$,$p\neq p^{\prime}$.
\end{lm}
{\it Proof}. By (\ref{cV-from-mnp}), we have
$$ n = c^2p,\quad n^{\prime} = c^2p^{\prime}.$$
If $p = p^{\prime}$, then $n = n^{\prime}$. Due to (\ref{Pell-equation}), $m^2 = (m^{\prime})^2$ and hence, 
$m = m^{\prime}$. If $p \neq p^{\prime}$, then $n \neq n^{\prime}$. According to (\ref{Pell-equation}) and (\ref{cV-from-mnp}), 
$$ m^2 - np =  m^2 - c^2p^2 = 1 = (m^{\prime})^2 - n^{\prime}p^{\prime} =  (m^{\prime})^2 - c^2(p^{\prime})^2.$$
Hence, $p \neq p^{\prime}$ implies  $m^2 \neq (m^{\prime})^2$ and $m \neq m^{\prime}$. $\Box$

\begin{lm}\label{unique-for-V}
Given an admissible $V$, there is a unique solution $(m,n,p)$ of equation (\ref{Pell-equation}).
\end{lm}
{\it Proof}. Due to (\ref{cV-from-mnp}), we have
$$ V = \frac{n}{m} = \frac{n^{\prime}}{m^{\prime}}. $$
By lemma \ref{coprime}, the pairs $(m,n)$, $(m^{\prime},n^{\prime})$ are coprime. Hence, $m = m^{\prime}$, 
$n = n^{\prime}$. It follows from (\ref{Pell-equation}) that $p = p^{\prime}$. $\Box$
\begin{cor}\label{fixed-V-correspond-one-c}
Given an admissible $V$, there is a unique admissible $c$ such that the pair $(c,V)$ is admissible.
\end{cor}

Corollary \ref{fixed-V-correspond-one-c} means that the parameter $c$ is a function of $V$, $c = c(V)$, that is, the parameter $V$ defines uniquely Lorentzian transformation of the torus. In contrast with lemma \ref{unique-for-V}, the following lemma means that a fixed admissible $c$ correspond infinitely many triples $(m,n,p)$.
\begin{lm}\label{fixed-c-correspond-infinitely-many-triples}
Given an admissible $c$, there are infinitely many different solutions $(m,n,p)$ of equation (\ref{Pell-equation}).
\end{lm}
{\it Proof}. Let us fix $c\stackrel{\rm def}{=}C_0$ corresponding to some solution $(M_0,N_0,P_0)$ of (\ref{Pell-equation}), $ C_0 = \sqrt{\frac{N_0}{P_0}} = \frac{\sqrt{M^2-1}}{P_0} $. We have to find infinitely many different solutions $(m,n,p)$ of (\ref{Pell-equation}) that satisfy the last equality with replacing $(M_0,N_0,P_0)$ by $(m,n,p)$.

Consider the equation
\begin{equation}\label{real-Pell-equation}
m^2 - ap^2 = 1,\quad \mbox{where}\quad a = \frac{N_0}{P_0}.
\end{equation}
Note that (\ref{real-Pell-equation}) in general is not a Pell's equation because $a$ is not necessary a natural number (see case 2) below).
Since we try to solve equation (\ref{Pell-equation}) with the condition $ \frac{n}{p} = \frac{N_0}{P_0} $, it is enough to prove that (\ref{real-Pell-equation}) has infinitely many different integer solutions $(m,p)$, $m\ge 2$, $p\in \mathbb{N}$. There are two cases: 1) $a$ is a natural number; 2) $a$ is not a natural number. In the case 1), 
$ \sqrt{a} = C_0 $ is irrational, according to lemma \ref{c-is-irrational}. Therefore, (\ref{real-Pell-equation}) is actually a Pell's equation (see, for example \cite{Sierpinski-book64}). Due to theorems 10.9.1 and 10.9.2 \cite{Sierpinski-book64}, a Pell's equation has infinitely many (nontrivial) positive solutions $(m_i,p_i)$. This gives the infinitely many different solutions $(m_i,n_i,p_i)$ of equation (\ref{Pell-equation}), where 
$ n_i = N_0\frac{p_i}{P_0} $.

In the case 2), replacing $a$ in (\ref{real-Pell-equation}) by $ a = \frac{N_0}{P_0}$, we get
$ m^2 - \frac{N_0}{P_0}p^2 = 1 $ or
\begin{equation}\label{another-real-Pell-equation}
m^2 - (N_0P_0)\left(\frac{p}{P_0}\right)^2 = 1.
\end{equation}
This is a Pell's equation with irrational 
$$ \sqrt{N_0P_0} = P_0\frac{N_0}{P_0} = C_0P_0$$
since $C_0$ is irrational (lemma \ref{c-is-irrational}). Hence, there are infinitely many (nontrivial) positive solutions $ (m_i,P_i) $, where $ P_i = \frac{p}{P_0} $. Again, this gives the infinitely many different solutions $(m_i,n_i,p_i)$ of equation (\ref{Pell-equation}), where $ p_i = P_iP_0 $, $ n_i = N_0P_i = N_0\frac{p_i}{P_0} $.
This completes the proof. $\Box$
\begin{cor}\label{spectra-is-infinite}
Every admissible $c$ has a countable spectra.
\end{cor}
\begin{lm}\label{points-of-accumulation}
Let $ V_1, V_2, \ldots , V_i, \ldots $ is the positive part of the spectra of $c$, $c = c(V_i)$. Then the set
$ \{V_1, V_2, \ldots , V_i, \ldots \} $ has a unique point of accumulation which is $c$.
\end{lm}
{\it Proof}. Since the set $ \{V_1, V_2, \ldots , V_i, \ldots \} $ is countable, it has at least one point of accumulation, say $ w\in [0;c] $. We have to show that $ w = c $. By lemma \ref{fixed-c-correspond-infinitely-many-triples}, there are infinitely many different solutions $(m,n,p)$ of equation (\ref{Pell-equation}) corresponding to $c$. Due to (\ref{Pell-equation}) and (\ref{cV-from-mnp}),
$$ \frac{V_i}{c} = \frac{\sqrt{m^2_i-1}}{m_i}\quad \mbox{for every}\quad i\ge 1. $$
If the sequence $ \{m_i\} $ is bounded, then there only finitely many $ V_i $. This contradicts to corollary \ref{spectra-is-infinite}. Then, for any unbounded subsequence $ m_{i_k}\to +\infty $, $ \frac{V_{i_k}}{c}\to 1 $ as 
$ k\to \infty $. Hence, $ V_{i_k}\to c $ as $ k\to \infty $. This proves the lemma. $\Box$
\begin{lm}\label{c-many}
Let $ C $ be a positive irrational number, which is the square root from some rational number, say $ \frac{N}{P} $. Then $ C $ is an admissible (light speed) parameter.
\end{lm}
{\it Proof}. By condition, $ C^2 = \frac{N}{P} $. Without loss of generality, we can assume that $ N $, 
$ P\in \mathbb{N} $. Consider the equation
\begin{equation}\label{again-Pell-equation}
m^2 - dr^2 = 1,\quad \mbox{where}\quad d = NP.
\end{equation}
Since $C$ is irrational, $CP = \sqrt{NP} = \sqrt{d}$ is also irrational. Hence, (\ref{again-Pell-equation}) is a Pell equation that has infinitely many positive (nontrivial) solutions $(m,r)$. Given the solution $(m_0,r_0)$, put
$n = r_0N$, $p = r_0P$. Then
$$ m_0^2 - np = m_0^2 - r_0^2NP = m_0^2 - dr_0^2 = 1. $$
Another words, $(m_0,n,p)$ is a positive (nontrivial) solution of equation (\ref{Pell-equation}). Since
$ C = \sqrt{\frac{N}{P}} = \sqrt{\frac{n}{p}} $, $ C $ is admissible. The lemma is proved. $\Box$
\begin{cor}
The set of admissible (light speed) parameters $c$ is dense on the real positeve halfline.
\end{cor}

{\it Proof of theorem \ref{arithmetical-properties}} follows from the lemmas \ref{c-is-irrational}, \ref{coprime} - \ref{c-many} above.

\section{Calculation of spectra}

Here, given admissible parameters (light speed) $c$, we present the algorithm calculating the corresponding spectra 
$ \ldots , -V_i, \ldots , -V_1, V_1, V_2, \ldots , V_i, \ldots $, where $c = c(V_i)$. Later on, $c$ if fixed. It is enough to present such algorithm for the positive part $ \{V_1, V_2, \ldots , V_i, \ldots \} $ of the spectra. Due to lemma \ref{points-of-accumulation}, $c$ is a unique point of accumulation for this positive part. Therefore, without loss of generality, we can assume that
$$ V_1 < V_2 < \ldots < V_i < V_{i+1} < \ldots . $$
According to (\ref{cV-from-mnp}) and (\ref{mnp-from-cV}), any admissible pair $(c, V_i)$ corespond uniquely the  triples $(m_i, n_i, p_i)$ of positive integer solutions of (\ref{Pell-equation}), such that
\begin{equation}\label{another-cV-from-mnp}
c = \sqrt{\frac{n_i}{p_i}},\quad V_i = \frac{n_i}{m_i},
\end{equation}
\begin{equation}\label{another-mnp-from-cV}
m_i = \frac{1}{\sqrt{1-\frac{V_i^2}{c^2}}},\quad n_i = \frac{V_i}{\sqrt{1-\frac{V_i^2}{c^2}}},\quad 
p_i = \frac{V_i}{c^2\sqrt{1-\frac{V_i^2}{c^2}}},
\end{equation} 
\begin{equation}\label{another-Pell-equation}
m_i^2 - n_ip_i = 1\quad m_i\ge 2, \quad n_i, p_i\in \mathbb{N}.
\end{equation}

In the set of arbitrary triples $(m, n, p)$, we introduce a partial order as follows. Say
$$ (m, n, p) < (m^{\prime}, n^{\prime}, p^{\prime})\quad \mbox{iff}
\quad m < m^{\prime}, n < n^{\prime}, p < p^{\prime}. $$
We also shall use a partial order for doubles:
$$ (m, n) < (m^{\prime}, n^{\prime})\quad \mbox{iff}
\quad m < m^{\prime}, n < n^{\prime}. $$
\begin{lm}\label{ordering}
Let $(m_i, n_i, p_i)$ corresponds to $(c, V_i)$. Then
$$ (m_1, n_1, p_1) < \ldots < (m_i, n_i, p_i) < (m_{i+1}, n_{i+1}, p_{i+1}) < \ldots . $$
\end{lm}
{\it Proof}. It follows from (\ref{another-mnp-from-cV}) and (\ref{another-Pell-equation}) that
$$ V_i^2 = \frac{m_i^2-1}{m_i^2}c^2 < V_{i+1}^2 = \frac{m_{i+1}^2-1}{m_{i+1}^2}c^2.$$
Hence, $ \frac{m_i^2-1}{m_i^2} <  \frac{m_{i+1}^2-1}{m_{i+1}^2}$, and $m_i < m_{i+1}$. Since
$ V_i = \frac{n_i}{m_i} < V_{i+1} = \frac{n_{i+1}}{m_{i+1}} $, $n_i < n _{i+1}$. Similarly,
$ c^2 = \frac{n_i}{p_i} = \frac{n_{i+1}}{p_{i+1}} $ implies $p_i < p_{i+1}$. $\Box$

Thus, to define the spectra for the given $c$ it is sufficient to find the first and minimal triple $ (m_1, n_1, p_1) $ and represent a recurrent formula in calculating $ (m_i, n_i, p_i) $, $i\ge 2$. Due to (\ref{another-cV-from-mnp}), each
$ (m_i, n_i, p_i) $ will correspond to admissible $V_i = \frac{n_i}{m_i}$ with $c = c(V_i) = \sqrt{\frac{n_i}{p_i}}$.

According to theorem \ref{arithmetical-properties}, $ c^2 = \frac{n_*}{p_*} $ for some $n_*$, $p_*\in \mathbb{N}$. Without loss of generality, we can assume that $n_*$, $p_*$ are coprime (otherwise, we divide $n_*$, $p_*$ on a common multiplier). Consider the equation
\begin{equation}\label{one-more-Pell-equation}
m^2 - dr^2 = 1,\quad \mbox{where}\quad d = c^2p^2_*.
\end{equation}
Since $c$ is irrational, $cp_* = \sqrt{c^2p^2_*} = \sqrt{d}$ is also irrational. Hence, (\ref{one-more-Pell-equation}) is a Pell equation. Due to \cite{Sierpinski-book64}, (\ref{one-more-Pell-equation}) has a minimal (in the set of natural numbers) solution, say $(m_0,r_0)$. Below we recall how one can get all solutions of Pell's equation and they have an order structure..

Put $ m_1 = m_0 $, $n_1 = r_0n_*$, $p_1 = r_0p_*$. Then
$$ m_1^2 - n_1p_1 = m_0^2 - r_0^2(n_*p_*) = m_0^2 - (c^2p_*^2)r_0^2 = m_0^2 - dr_0^2 = 1. $$
Moreover, $ c^2 = \frac{n_*}{p_*} = \frac{n_1}{p_1} $. Hence, the triple $ (m_1, n_1, p_1) $ correspond to some term in the spectra of $c$. Let us show that $ (m_1, n_1, p_1) $ is a minimal triple.

Suppose the contrary. Then there is the triple $ (m^{\prime},n^{\prime},p^{\prime}) < (m_1, n_1, p_1)$ with
$ c^2 = \frac{n^{\prime}}{p^{\prime}} $, $ (m^{\prime})^2 - n^{\prime}p^{\prime} = 1 $. We have
$ \frac{n^{\prime}}{p^{\prime}} = \frac{n_*}{p_*} $, and so $ n^{\prime} = \frac{p^{\prime}}{p_*}n_* $. Since $n_*$ and $p_*$ are coprime, $ \frac{p^{\prime}}{p_*} $ is a natural number, say $ r^{\prime} $. The straight calculation shows that $ (m^{\prime}, r^{\prime}) $ is a solution of (\ref{one-more-Pell-equation}):
$$ (m^{\prime})^2 - (c^2p^2_*)(r^{\prime})^2 = (m^{\prime})^2 - (c^2p^2_*)(\frac{p^{\prime}}{p_*})^2 =
(m^{\prime})^2 - c^2(p^{\prime})^2  = (m^{\prime})^2 - n^{\prime}p^{\prime} = 1.$$
Since $ (m_1, p_1) $ is the minimal solution of (\ref{one-more-Pell-equation}), $ m_1 < m^{\prime} $. This contradicts to inequality $ (m^{\prime},n^{\prime},p^{\prime}) < (m_1, n_1, p_1)$ that means, in particular, 
$ m_1 > m^{\prime} $.

{\bf Resume}. To get the minimal triple $ (m_1, n_1, p_1) $, one need to take the minimal (in the set of natural numbers) solution $(m_0,r_0)$ of (\ref{one-more-Pell-equation}). Then
$$ m_1 = m_0 ,\quad n_1 = r_0n_*,\quad p_1 = r_0p_*$$
gives the minimal triple $ (m_1, n_1, p_1) $. This triple corresponds to the minimal term $V_1$ of the spectra of $c$ by (\ref{another-cV-from-mnp}).

Remark, that if $ c^2 = \frac{n_*}{p_*} $ with  coprime $n_*$, $p_*$ then, in general, $n_1\neq n_*$ and 
$p_1\neq p_*$. For example, if $ c^2 = \frac{1}{2} $, then $ (m_1, n_1, p_1) = (3,2,4)$.

Let us recall the solution of Pell's equation (\ref{one-more-Pell-equation}). Due to theorem 3.5.7 \cite{Sierpinski-book64}, if the period of the continued fraction of $d$ consists of an even numbers $s$ of terms, then the numerator and the denominator of the $(ns-1)$-th convergent, $n\in \mathbb{N}$, form a solution of (\ref{one-more-Pell-equation}). Moreover, all the solutions are obtained in this way. From this we see that the solution in the least natural numbers is given by the $(s-1)$-th convergent. Since numerators and denominators form increasing sequences, the solutions of (\ref{one-more-Pell-equation}) are endowed with the natural order:
$$ (P_{s-1},Q_{s-1}) <  (P_{2s-1},Q_{2s-1}) < \ldots (P_{ns-1},Q_{ns-1}) < \ldots .$$
Due to theorem 3.5.8 \cite{Sierpinski-book64}, if the period of the continued fraction of $d$ consists of an odd numbers $s$ of terms, then the numerator and the denominator of the $(2ns-1)$-th convergent, $n\in \mathbb{N}$, form a solution of (\ref{one-more-Pell-equation}). In this case, all the solutions are also obtained in this way. From this we see that the solution in the least natural numbers is given by the $(2s-1)$-th convergent. Again, the solutions of (\ref{one-more-Pell-equation}) are endowed with the natural order:
$$ (P_{2s-1},Q_{2s-1}) <  (P_{4s-1},Q_{4s-1}) < \ldots (P_{2ns-1},Q_{2ns-1}) < \ldots .$$

It is easy to see that the system of equations $ m^2 - np = 1 $, $ c^2 = \frac{n}{p} $ is equivalent to the system
 $ m^2 - c^2p^2 = 1 $, $ n = c^2p $. Therefore, we consider the equation
\begin{equation}\label{quasi-Pell}
X^2 - c^2Y^2 = 1.
\end{equation}
Let us introduce an operation $\bigotimes$ on pairs of numbers as follows:
$$ (u,v)\bigotimes (x,y) = (ux + c^2vy, vx + uy). $$
\begin{lm}\label{structure-of-solutions}
Let $ (u,v)$ and $(x,y) $ are solutions (not necessary, integer) of (\ref{quasi-Pell}). Then 
$ (u,v)\bigotimes (x,y) = (ux + c^2vy, vx + uy) $ is the solution of (\ref{quasi-Pell}) as well.
\end{lm}
{\it Proof}. The straight calculation gives
$ (ux + c^2vy)^2 - c^2(vx + uy)^2 = (x^2 - c^2y^2)(u^2 - c^2v^2) = 1$.
$\Box$
\begin{cor}\label{1-from-structure-of-solutions}
Let $ (m,p)$ and $(m^{\prime},p^{\prime}) $ are integer solutions (not necessary, positive) of (\ref{quasi-Pell}). Then $(m,p)\bigotimes (m^{\prime},p^{\prime})$ is the solution of (\ref{quasi-Pell}) as well.
\end{cor}
Taking in mind that the minimal solution  $(m_1,p_1)$ of equation (\ref{quasi-Pell}) consists of positive $m_1$ and $p_1$, we get
\begin{cor}\label{2-from-structure-of-solutions}
Let $ (m_1,p_1)$ be the minimal solution of (\ref{quasi-Pell}). Then 
$$\underbrace{(m_1,p_1)\bigotimes \cdots \bigotimes (m_1,p_1)}_{k \mbox {times }}
\stackrel{\rm def}{=}(m_1,p_1)^k$$
 is the solution of (\ref{quasi-Pell}) as well. Moreover
$$ (m_1,p_1) < (m_1,p_1)^2 < \ldots (m_1,p_1)^k < \ldots . $$
\end{cor}

Consider the transformation $G: \mathbb{R}^2 \to \mathbb{R}^2$ of the type
$$ x^{\prime} = m_1x - c^2p_1y, \quad y^{\prime} = -p_1x + m_1y, $$
where $ (m_1,p_1)$ is the minimal solution of (\ref{quasi-Pell}). The determinant of $G$ equals 1, so $G$ is an orientation preserving diffeomorphism of $\mathbb{R}^2$. Moreover, due to lemma \ref{structure-of-solutions}, the hyperbolas $x^2 - c^2y^2 = 1$ denoted by $H$ is invariant under $G$.

Let $M_k$ be the point with the coordinates $ (m_1,p_1)^k $. According to corollary \ref{2-from-structure-of-solutions}, $M_k\in H$ for every $k\in \mathbb{N}$. The point $M_0$ correspond to the trivial solution $(1,0)$ of $(m^{\prime},p^{\prime})$, and is the vertex of right side branch of the hyperbolas $H$.
\begin{lm}\label{all-solutions}
Let $(m^{\prime},p^{\prime})$ be a nontrivial positive integer solution of (\ref{quasi-Pell}). Then $(m^{\prime},p^{\prime})$ is $ (m_1,p_1)^k $ for some $k\in \mathbb{N}$.
\end{lm}
{\it Proof}. Suppose the contrary. Then $ (m_1,p_1)^k < (m^{\prime},p^{\prime}) < (m_1,p_1)^{k+1}$ for some $k\in \mathbb{N}$. This means that the point $M^{\prime}$ with the coordinates $(m^{\prime},p^{\prime})$ is between points $M_k$, $M_{k+1}$ on $H$. Recall that $H$ is invariant under $G$ and, hence, under $G^{-k}$. Under the diffeomorphism $G^{-k}$, the points $M_k$, $M_{k+1}$ are mapped into the points $M_0$ and $M_1$ respectively. Since $G$ preserves orientation, $M^{\prime}$ is mapped to the point $G^{-k}(M^{\prime})$ between $M_0$ and $M_1$. Due to corollary \ref{1-from-structure-of-solutions}, $G^{-k}(M^{\prime})$ corresponds to an integer solution of (\ref{quasi-Pell}). Since $G^{-k}(M^{\prime})$ is between $M_0$ and $M_1$, this solution is nontrivial, positive, and less that $(m_1,p_1)$. This contradiction proves the lemma. $\Box$

As a consequence of lemma \ref{all-solutions}, we get that the triples $ (m_k, n_k, p_k) $ are obtained from the following relations:
$$ (m_k, p_k) = (m_1,p_1)^k, \quad n_k = \frac{m_k^2-1}{p_k}, \quad k\in \mathbb{N}. $$
Let us get direct recurrent formulas.
\begin{lm}\label{recurrent-formulas}
$ (m_k, p_k) = (m_1,p_1)^k $ iff 
\begin{equation}\label{RecurrentFormulas}
 m_k = m_1m_{k-1} + c^2p_1p_{k-1}, \quad  p_k = m_1p_{k-1} + p_1m_{k-1}, \quad k\ge 1, 
\end{equation}
where $m_0 = 1$, $p_0 = 0$.
\end{lm}
{\it Proof}. Let $ (m_k, p_k) = (m_1,p_1)^k $. Obviously, (\ref{RecurrentFormulas}) takes place  for $k = 1$. Suppose (\ref{RecurrentFormulas}) are true for $1$, $\ldots$, $k-1 \ge 1$, and show that (\ref{RecurrentFormulas}) is true for $k$. By induction suggestion, we have
$$ m_{k-1} = m_1m_{k-2} + c^2p_1p_{k-2}, \quad  p_{k-1} = m_1p_{k-2} + p_1m_{k-2}. $$
Since 
$ (m_k, p_k) = (m_1,p_1)^k = (m_1,p_1)\bigotimes (m_1,p_1)^{k-1} = (m_1,p_1)\bigotimes (m_{k-1},p_{k-1}) $,
(\ref{RecurrentFormulas}) follows by definition of $\bigotimes$. The converse statement is proved similarly.
$\Box$

{\bf Resume.} After finding the minimal triple $ (m_1,n_1p_1) $, the recurrent formulas
$$ m_k = m_1m_{k-1} + c^2p_1p_{k-1}, \quad  p_k = m_1p_{k-1} + p_1m_{k-1}, \quad 
n_k = \frac{m_k^2-1}{p_k} $$
give the $k$-th triple $ (m_k, n_k, p_k) $. This triple corresponds to the $k$-th term $(V_{-k}, V_k) = (-V_k, V_k)$ of the spectra of $c$ by (\ref{another-cV-from-mnp}): $ V_k = \frac{n_k}{m_k} $.

Remark that the relation $ (m_k, p_k) = (m_1,p_1)^k $ is equivalent to 
$$ m_k + cp_k = (m_1 + cp_1)^k .$$
This can be proved similarly to the proof of lemma \ref{recurrent-formulas}.

\bigskip

2875 COWLEY WAY (1015), SAN DIEGO, CA 92110, USA

{\it E-mail address}: saranson@yahoo.com
\medskip

DEPARTMENT OF APPL. MATHEMATICS, NIZHNY NOVGOROD STATE TECHNICAL UNIVERSITY, NIZHNY NOVGOROD, RUSSIA

{\it E-mail address}: zhuzhoma@mail.ru

{\it Current e-mail address}: zhuzhoma@maths.univ-rennes1.fr


\begin{thebibliography}{99}

\bibitem{AdamsStuch97}
{\bf Adams S., Stuck G.} The isometry groups of Lorentz manifolds I, II. {\it Invent. Math.}, {\bf 129}(1997), 239-261; 263-287.

\bibitem{BeemEhrlich-book81}
{\bf Beem J., Ehrlich P.} {\it Global Lorentzian Geometry}. Monographs and Textbooks in Pure and Appl. Math., {\bf 67}(1981).

\bibitem{DAmbra88}
{\bf D'Ambra G.} Isometry groups of Lorentz manifolds I, II. {\it Invent. Math.}, {\bf 95}(1988), 555-565.

\bibitem{Robinson-book99}
{\bf Robinson C.} {\it Dynamical Systems, Stability, Symbolic Dynamics, and Chaos}, second edition. CRC Press, {\bf 1999}.

\bibitem{Sierpinski-book64}
{\bf Sierpinski W.} {\it Elementary Theory of Numbers}. Warsaw, {\bf 1964}.

\bibitem{Smale67}
{\bf Smale S.} {\it Differentiable dynamical systems.} Bull. Amer. Math. Soc.,
1967, 73, 1, 741-817.

\bibitem{Wolf-book72}
{\bf Wolf J.} {\it Spaces of Constant Curvature}. Univ. of Carolina, Berkley, {\bf 1972}.



\end{thebibliography}
\end{document}